\def\bbbr{{\rm I\!R}} 

\def\bbbn{{\rm I\!N}} 
\def\bbbh{{\rm I\!H}}

\def\bbbone{{\mathchoice {\rm 1\mskip-4mu l} {\rm 1\mskip-4mu l}
{\rm 1\mskip-4.5mu l} {\rm 1\mskip-5mu l}}}

\def\bbbc{{\mathchoice {\setbox0=\hbox{$\displaystyle\rm C$}\hbox{\hbox
to0pt{\kern0.4\wd0\vrule height0.9\ht0\hss}\box0}}
{\setbox0=\hbox{$\textstyle\rm C$}\hbox{\hbox
to0pt{\kern0.4\wd0\vrule height0.9\ht0\hss}\box0}}
{\setbox0=\hbox{$\scriptstyle\rm C$}\hbox{\hbox
to0pt{\kern0.4\wd0\vrule height0.9\ht0\hss}\box0}}
{\setbox0=\hbox{$\scriptscriptstyle\rm C$}\hbox{\hbox
to0pt{\kern0.4\wd0\vrule height0.9\ht0\hss}\box0}}}}

\def\bbbe{{\mathchoice {\setbox0=\hbox{\smalletextfont e}\hbox{\raise
0.1\ht0\hbox to0pt{\kern0.4\wd0\vrule width0.3pt
height0.7\ht0\hss}\box0}}
{\setbox0=\hbox{\smalletextfont e}\hbox{\raise
0.1\ht0\hbox to0pt{\kern0.4\wd0\vrule width0.3pt
height0.7\ht0\hss}\box0}}
{\setbox0=\hbox{\smallescriptfont e}\hbox{\raise
0.1\ht0\hbox to0pt{\kern0.5\wd0\vrule width0.2pt
height0.7\ht0\hss}\box0}}
{\setbox0=\hbox{\smallescriptscriptfont e}\hbox{\raise
0.1\ht0\hbox to0pt{\kern0.4\wd0\vrule width0.2pt
height0.7\ht0\hss}\box0}}}}

\def\bbbq{{\mathchoice {\setbox0=\hbox{$\displaystyle\rm Q$}\hbox{\raise
0.15\ht0\hbox to0pt{\kern0.4\wd0\vrule height0.8\ht0\hss}\box0}}
{\setbox0=\hbox{$\textstyle\rm Q$}\hbox{\raise
0.15\ht0\hbox to0pt{\kern0.4\wd0\vrule height0.8\ht0\hss}\box0}}
{\setbox0=\hbox{$\scriptstyle\rm Q$}\hbox{\raise
0.15\ht0\hbox to0pt{\kern0.4\wd0\vrule height0.7\ht0\hss}\box0}}
{\setbox0=\hbox{$\scriptscriptstyle\rm Q$}\hbox{\raise
0.15\ht0\hbox to0pt{\kern0.4\wd0\vrule height0.7\ht0\hss}\box0}}}}

\def\bbbt{{\mathchoice {\setbox0=\hbox{$\displaystyle\rm
T$}\hbox{\hbox to0pt{\kern0.3\wd0\vrule height0.9\ht0\hss}\box0}}
{\setbox0=\hbox{$\textstyle\rm T$}\hbox{\hbox
to0pt{\kern0.3\wd0\vrule height0.9\ht0\hss}\box0}}
{\setbox0=\hbox{$\scriptstyle\rm T$}\hbox{\hbox
to0pt{\kern0.3\wd0\vrule height0.9\ht0\hss}\box0}}
{\setbox0=\hbox{$\scriptscriptstyle\rm T$}\hbox{\hbox
to0pt{\kern0.3\wd0\vrule height0.9\ht0\hss}\box0}}}}

\def\bbbs{{\mathchoice
{\setbox0=\hbox{$\displaystyle     \rm S$}\hbox{\raise0.5\ht0\hbox
to0pt{\kern0.35\wd0\vrule height0.45\ht0\hss}\hbox
to0pt{\kern0.55\wd0\vrule height0.5\ht0\hss}\box0}}
{\setbox0=\hbox{$\textstyle        \rm S$}\hbox{\raise0.5\ht0\hbox
to0pt{\kern0.35\wd0\vrule height0.45\ht0\hss}\hbox
to0pt{\kern0.55\wd0\vrule height0.5\ht0\hss}\box0}}
{\setbox0=\hbox{$\scriptstyle      \rm S$}\hbox{\raise0.5\ht0\hbox
to0pt{\kern0.35\wd0\vrule height0.45\ht0\hss}\raise0.05\ht0\hbox
to0pt{\kern0.5\wd0\vrule height0.45\ht0\hss}\box0}}
{\setbox0=\hbox{$\scriptscriptstyle\rm S$}\hbox{\raise0.5\ht0\hbox
to0pt{\kern0.4\wd0\vrule height0.45\ht0\hss}\raise0.05\ht0\hbox
to0pt{\kern0.55\wd0\vrule height0.45\ht0\hss}\box0}}}}

\def\sq{\hbox{\rlap{$\sqcap$}$\sqcup$}}

\def\qed{\ifmmode\sq\else{\unskip\nobreak\hfil
\penalty50\hskip1em\null\nobreak\hfil\sq
\parfillskip=0pt\finalhyphendemerits=0\endgraf}\fi}

\def\Proof{\removelastskip\vskip\baselineskip\noindent{\it
Proof.\quad}\ignorespaces}

\newtheorem{exer}{Exercise}

\newtheorem{theor}{Theorem}
\def\theorem#1#2{
		\begin{theor}\label{#1}
		#2
		\end{theor}
}

\newtheorem{propos}{Proposition}

\newtheorem{defini}{Definition}
\def\definition#1#2{
		\begin{defini}\label{#1}
		#2
		\end{defini}
}

\newtheorem{lemm}{Lemma}
\def\lemma#1#2{
		\begin{lemm}\label{#1}
		#2
		\end{lemm}
}

\newtheorem{corol}{Corollary}
\def\corollary#1#2{
		\begin{corol}\label{#1}
		#2
		\end{corol}
}

\def\A{A}
\def\A{A}   
\def\B{B}
\def\CC{C}
\def\C{{\bbbc}}
\def\Dilat#1{D_{#1}} 

\def\F{F}

\def\Halfplane{\bbbh} 

\def\H{H}
\def\H{H} 

\def\Id{\bbbone} 
\def\InvWavetrans{{\cal M}}

\def\K{K} 
\def\L#1{L^{#1}}   

\def\M{M}
\def\Ndim{N}

\def\N{{\bbbn}}  
\def\O{O}  

\def\Pol{P}

\def\Rplus{{\R_+}}
\def\R{{\bbbr}} 

\def\Schwarz{S}  

\def\Translat#1{T_{#1}} 
\def\T{T}
\def\T{T}
\def\U{U}
\def\M{M}

\def\Wavetrans{{\cal W}} 

\def\x{x}

\def\aa{{\a^\prime}} 
\def\aaa{{a^{\prime\prime}}}
\def\abs#1{\left\vert #1 \right\vert} 

\def\a{a}       
\def\bbb{{b^{\prime\prime}}}
\def\bb{{\b^\prime}}

\def\b{b}      
\def\conj#1{\overline{#1}} 
\def\const{c}

\def\d{d}    

\def\e{e}   
\def\fou#1{\widehat{#1}} 
\def\f{f}

\def\g{g}    
\def\gg{{g^\prime}}
\def\h{h}    
\def\hh{{\h^\prime}}

\def\intinf{\int_{-\infty}^{+\infty}} 
\def\inv#1{{1\over #1}}

\def\i{i} 

\def\k{k}   

\def\l{l}    
\def\max{{\rm max}\,}

\def\m{m}    

\def\norm#1{\left\Vert #1 \right\Vert}
\def\n{n}    
\def\of#1{\left( #1 \right)}
\def\o{o}    

\def\p{p}
\def\q{q}

\def\r{r}

\def\scalarproduct#1#2{\left\langle#1\mid#2\right\rangle}

\def\suchthat{:}

\def\supp{\hbox{\rm supp}\,}
\def\s{s}

\def\t{t}      
\def\u{u}
\def\x{x}

\def\y{y}

\font\twelve=cmbx10 at 15pt
\font\ten=cmbx10 at 12pt

\documentstyle[12pt,epsf]{article}

\def\dist{\hbox{dist}}
\begin{document}

\begin{titlepage}

\begin{center}

\renewcommand{\thefootnote}{\fnsymbol{footnote}}

{\ten Centre de Physique Th\'eorique\footnote{
Unit\'e Propre de Recherche 7061} - CNRS - Luminy, Case 907}
{\ten F-13288 Marseille Cedex 9 - France }

\vspace{1 cm}

{\twelve Some directional microlocal classes \\
defined using  \\
wavelet transforms}

\vspace{0.3 cm}

\setcounter{footnote}{0}
\renewcommand{\thefootnote}{\arabic{footnote}}

{\bf Matthias HOLSCHNEIDER\footnote{e-mail: hols@marcptnx1.univ-mrs.fr}}

\vspace{2 cm}

{\bf Abstract}

\end{center}
In this short paper we discuss  how the position - scale
half-space ofwavelet analysis may be cut into different regions. We
discuss conditions
under which they are independent in the sense that the T\"oplitz
operators
associated with their characteristic functions commute modulo smoothing
operators. This shall be used to define microlocal classes of
distributions
having a well defined behavior along  lines in wavelet space. This
allows us the description of singular and
regular directions in distributions. As an application we discuss
elliptic regularity for these microlocal classes for domains with
cusp-like singularities.
\vspace{2 cm}

\bigskip

\noindent Number of figures : 0

\bigskip

\noindent October 1995

\noindent CPT-95/P3248

\bigskip

\noindent anonymous ftp or gopher: cpt.univ-mrs.fr

\end{titlepage}

\section{Introduction}

The classical definition of local singular directions of a distribution
 $\eta$ is given by the wavefront set (e.g. \cite{HORMANDERI}). At
a given point
 $\x$ it is roughly speaking the cone of all directions in which the
 Fourier transform of the localized distribution $\phi\cdot\eta$ does
 not decay rapidly where $\phi$ is any smooth function that is supported
 by some neighborhood of $\x$. More precisely for fixed
 $\phi\in\CC^\infty_0(\R^\n)$, $\x\in\supp\phi$, a direction
 $\xi\in\R^\n-\{0\}$ is regular if in some conic neighborhood
 $\gamma\owns\xi$ we have
$$
\k\in\gamma\Rightarrow
\abs{(\phi\cdot\eta)\fou{\ }(\k)} \leq \const_\Ndim (1+\abs\k)^{-\Ndim},
\Ndim = 0,1,2,\dots,
$$
The complement of the regular directions is denoted by
$\Sigma_{\x,\phi}$.
 The singular directions at $\x$ are then defined as
$$
\Sigma_\x = \cap_{\phi}\Sigma_{\x,\phi},
$$
where the the intersection is over all $\phi\in\CC^\infty_0(\R^\n)$ with
 $\x\in\supp\phi$.
However this concept  of singular and regular directions does not
always fit with what one would intuitively call a regular direction
in a distribution.
 What we mean is best illustrated by an example in $\R^2$: consider
the set
 $\K=\{(\x,\y)\in\R^2 \suchthat \x\geq0, \abs\y\leq\x^2\}$ and let
$\chi$ be a function that is of very low regularity in the
complement of the cusp $\K$, whereas inside the cusp it is smooth.
It is plain to see that there is no direction in which the Fourier
 transform of the localized function $\phi\chi$ decays rapidly and so all
 directions are singular. This however is in contradiction with our
 intuition in the sense that if we approach the singularity along any
 path contained in the  set $\K$ no irregularity
 is to be noticed and one would like to call the  direction of the
cusp regular. Let us refine this example a little more and let us
consider a functions that inside $\K$ looks like
$$
\sin(\x^{-\alpha}),\quad\alpha>1
$$
whereas outside it is very irregular. This function although smooth
inside $\K$ becomes less and less smooth as we approach the top of
the cusp due to stronger and stronger oscillations. We therefore
see, that  the directional smoothness we want to study should have
at least one more parameter which corresponds to the way how the
oscillations  accumulate. Therefore let us look at details of size
$\a$ inside the cusp at distance $\b$ from the top. In order to see
the amount of details increase---that is in order to recover the
degree of increasing non-smoothnes near the top of the cusp---we
must choose $\b$ of the order of magnitude of $\a^{\alpha}$,
Otherwise we look at to small a scale compared to the distance, and
at this scale our function seems to be uniformly smooth locally.

This is clearly only a very vague statement. To give a more precise
definition we have to introduce the wavelet transform.
We shall be very brief and we refer to the literature for  a more
detailed
discussion (e.g. \cite{DAUBECHIESI}, \cite{MEYERIAII},
\cite{COIFMANMEYER}, \cite{HOLSBOOKWAVE}).

Let $\g\in\Schwarz(\R^\n)$, the class of Schwarz of rapidly decaying
 functions.  In addition suppose that $\g$ has all moments vanishing
$$
\int\d\x\,\x^\m\g(\x) =0
$$
for all multi indices $\m$.
Then the wavelet transform of $\s\in\L\p(\R^\n)$ with respect to $\g$ is
 defined as the following convolutions
\begin{equation}\label{eq: singdirect@waveletdef}
\Wavetrans_\g\s(\b,\a) = \Wavetrans[\g,\s](\b,\a) =
(\tilde\g_\a\ast\s)(\b)   = \int\d\x\inv{\a^\n}\conj\g\of{\x-\b\over\a}\,
\s(\x).
\end{equation}
with $\a>0$ and $\b\in\R^\n$.

Here we have introduced the following notations, that we shall use
in the sequel
$$
\tilde\g(\t)=\conj\g(-\t),\  \g_\a=\g(\cdot/\a)/\a^\n,
\g_{\b,\a}=\g_\a(\cdot-\b).
$$
The wavelet transform thus maps functions over the real line to
functions over the open half-space $\Halfplane^\n=\{(\b,\a)
\suchthat \b\in\R^\n, \a>0\}$.

 From the definition it is clear that the wavelet transform is a sort of
 mathematical microscope whose position is fixed by $\b$ and whose
 enlargement is given by $1/\a$ or to put it differently
 $\Wavetrans[\g,\s](\b,\a)$ is
 obtained by \lq\lq looking at $\s$ at position $\b$ and at scale
$\a$\rq\rq.
As a general statement one can say that
 local regularity of $\s$ is mirrored in a certain speed of decay of
 $\Wavetrans_\g\s$. So is for instance a uniform (in $\b$) decay of
 $\O(\a^\infty)$ as $\a\to0$, of the wavelet coefficients
equivalent to $\CC^\infty$
 regularity of $\s$. More quantitative information is available. So
a uniform decay of $\O(\a^\alpha)$ with $\alpha\in(0,1)$ is
equivalent to $\s\in\Lambda^\alpha$, the space of H\"older
continuous functions of exponenet $\alpha$.

A temptative definition of regular  direction at $\x$  is therefore any
 direction $\xi$ for which the wavelet transform decays faster than any
 power of $\a$ if the microscope approaches the singularity along a path
 that is tangent to  $\xi$ in $\x$ in such a way, that it looks at
a scale that is small compared to the distance to $\x$. That is we
say, vaguely speaking, a direction is regular if along a parabolic
line we have rapid decay of the wavelet coefficients
$$
\Wavetrans_\g\s(\lambda\xi,\lambda^\gamma) \leq
\O(\lambda^\infty),\quad (\lambda\to0).
$$
This idea will  be made more
 precise in section \ref{sec: singdirect@wavfst}.  In particular
the definition will we modified in such a way that it becomes
independent of the choice of the wavelet $\g$.

\section{The basic formulas of continuous wavelet analysis}
For the convenience of the reader we shall list here the basic
formulas of wavelet analysis. We limit ourselves to formal
expressions. They actually have a precise meaning when we consider
the wavelet analysis in $\Schwarz_0(\R^\n)$ or
$\Schwarz_0^\prime(\R^\n)$ (see below).

Let $\f$ be a complex valued function over $\R^\n$. Let $\g$ be
another such function. The wavelet transform of $\f$ with respect to
the analyzing wavelet $\g$ is defined through (we write $\d\x$ for
$\n$-dimensional Lebesgue measure)
$$
\Wavetrans_\g\f(\b,\a) =
\int_{\R^\n}\d\x\,\inv{\a^\n}\conj\g\of{\x-\b\over\a}\, \f(\x),\quad
\b\in\R^\n, \a>0.
$$
We also write $\Wavetrans[\g,\f](\b,\a)$ instead of
$\Wavetrans_\g\f(\b,\a)$.
Here $\b\in\R^\n$ is a position parameter and $\a\in\Rplus$ is a
scale parameter. The wavelet transform of a function over $\R^\n$ is
thus a function over the position scale half-space
$\Halfplane^\n=\R^\n\times\Rplus$.

If we introduce the dilation ($\Dilat\a$) and translation operators
($\Translat\b$)
$$
\Translat\b\s(\x) = \s(\x-\b),\quad
\Dilat\a\s(\x) = \s(\x/\a)/\a^\n,
$$
then we may also write the wavelet transform as a family of scalar
products
$$
\Wavetrans_\g\f(\b,\a) = \scalarproduct{\g_{\b,\a}}{\f},\quad
\g_{\b,\a}= \Translat\b\Dilat\a\g,
$$
or as a family of convolutions indexed by a scale parameter
$$
\Wavetrans_\g\f(\b,\a) = (\tilde\g_\a\ast\f)(\b),\quad
\tilde\g_\a=\Dilat\a\tilde\g,\quad\tilde\g(\x) = \conj\g(-\x).
$$
The convolution product is defined as usual
$$
(\s\ast\r)(\x) = \int_{\R^\n}\d\y \s(\x-\y)\r(\y) = (\r\ast\s)(\x)
$$

If we introduce the Fourier transform
$$
\fou\s(\k) = \int_{\R^\n}\d\x\e^{-\i\k\x}\s(\x),
\quad \s(\x)  = {1\over (2\pi)^\n}\int_{\R^\n}\d\k\e^{\i\k\x}\fou\s(\k)
$$
then the wavelet transform may also be written as
$$
\Wavetrans_\g\f(\b,\a) =
{1\over(2\pi)^\n}\int_{\R^\n}\d\k\,\conj{\fou\g(\a\k)}\,\e^{\i\b\k}\,
\fou\f(\k).
$$

The wavelet synthesis $\InvWavetrans$ maps functions over the
position - scale half-space to functions over $\R^\n$. Let
$\r=\r(\b,\a)$ be a complex valued function over $\Halfplane^\n$ and
$\h$ a function over $\R^\n$. Then the wavelet synthesis of $\r$
with respect to the synthesizing wavelet $\h$ is defined as
$$
\InvWavetrans_\h\r(\x) = \int_{\Halfplane^\n}{\d\b\d\a\over\a}\,
\r(\b,\a)\,
\inv{\a^\n}\, \h\of{\x-\b\over\a}.
$$

\subsection{Relation between $\Wavetrans$ and $\InvWavetrans$}

We now list some relations between $\Wavetrans$ and
$\InvWavetrans$.The wavelet synthesis is the adjoint of the wavelet
transform---both with respect to the same
wavelet---$\Wavetrans_\g^\ast=\InvWavetrans_\g$
$$
\int_{\Halfplane^\n}{\d\b\d\a\over\a} \conj{\Wavetrans_\g\s(\b,\a)}
\r(\b,\a) =
\int_{\R^\n}\d\x\, \conj{\s(\x)}\, \InvWavetrans_\g\r(\x).
$$
The combination $\InvWavetrans_\h\Wavetrans_\g$ reads in Fourier space
$$
\InvWavetrans_\h\Wavetrans_\g : \fou\s(\k) \mapsto  \m_{\g,\h}(\k)
\fou\s(\k),\quad
\m_{\g,\h}(\k) = \int_{0}^\infty
{\d\a\over\a} \conj{\fou\g(\a\k)}\, \fou\h(\a\k)
$$
Note that the Fourier multiplier $\m_{\g,\h}$ depends only on the
direction of $\k$, $\m_{\g,\h}=\m_{\g,\h}(\k/\abs\k)$. This is
because the measure $\d\a/\a$ is scaling invariant.
In case that $\g$ and $\h$ are such that
$$
\m_{\g,\h}(\k) = \int_{0}^\infty
{\d\a\over\a} \conj{\fou\g(\a\k)}\, \fou\h(\a\k) = \const_{\g,\h}
$$ with $0<\abs{\const_{\g,\h}} < \infty$,  we say that $\g$, $\h$
are an analysis reconstruction pair, or that $\h$ is a
reconstruction wavelet for $\g$. In this case we have
$$
\InvWavetrans_\h\Wavetrans_\g = \const_{\g,\h} \Id.
$$
We say that $\g$ is strictly  admissible if $\g$ is its own
reconstruction wavelet, or (what is the same) if
$$
\forall\k\in\R^\n\backslash\{0\}:\quad\int_{0}^\infty
{\d\a\over\a} \abs{\fou\g(\a\k)}^2 =  \const_\g
$$

A sufficient condition for $\g$ to have a reconstruction wavelet
$\r$ is that for some $\const>1$ we have
$$
\const^{-1} \leq \int_{0}^\infty
{\d\a\over\a} \abs{\fou\g(\a\k)}^2 \leq \const.
$$
If this condition holds, we call $\g$ admissible. In this case the
following function $\r$ will be a reconstruction wavelet for $\g$
\begin{equation}\label{recwaveform}
\fou\r(\k) = \fou\g(\k) \Big / \sqrt{\int_{0}^\infty
{\d\a\over\a} \abs{\fou\g(\a\k)}^2}.
\end{equation}

If $\g$ and $\h$ are an analysis reconstruction pair, then the
following formula holds
$$
\int_{\Halfplane^\n}{\d\b\d\a\over\a} \conj{\Wavetrans_\g\s(\b,\a)}
\Wavetrans_\h\r(\b,\a) = \const_{\g,\h}
\int_{\R^\n}\d\x\, \conj\s(\x)\,  \r(\x).
$$
In particular if $\g$ is strictly admissible, then we have
conservation of energy ($\const_{\g} = \const_{\g,\g}$)
$$
\int_{\Halfplane^\n}{\d\b\d\a\over\a}
\abs{\Wavetrans_\g\s(\b,\a)}^2  = \const_{\g}
\int_{\R^\n}\d\x\, \abs{\s(\x)}^2.
$$

 Let $\Pi$ and $\r$ be functions over the half-space
$\Halfplane^\n$.We define a (non-commutative) convolution for
functions
  over the half-space via
\begin{equation}\label{convdefhlfplsne}
(\Pi\ast\r)(\b,\a) = \int_{\Halfplane^\n}{\d\bb\d\aa\over\aa}\,
 {1\over\aa^\n}\Pi\of{{\b-\bb\over\aa},{\a\over\aa}}\, \r(\bb,\aa).
\end{equation}

Formally we can write
$$
(\Wavetrans_\g\InvWavetrans_\h\s)(\b,\a) = \Pi_{\g,\h}\ast\s(\b,\a),\quad
\Pi_{\g,\h} = \Wavetrans_\g\h.
$$

In the case where $\g$, and $\h$ are an analysis reconstruction
pair with $\const_{\g,\h}=1$ we clearly have
$(\Wavetrans_\g\InvWavetrans_\h)^2 = \Wavetrans_\g\InvWavetrans_\h$
and hence
$$
\Pi_{\g,\h}\ast\Pi_{\g,\h} = \Pi_{\g,\h}.
$$
The mapping $\r\mapsto \Pi_{\g,\h}\ast\r$ is a projector into the
range of the wavelet transform $\Wavetrans_\g$. In case that $\g$ is
strictly admissible with $\const_{g,\g}=1$,  we have that
$\Pi_{\g,\g}$ is an orthogonal projector.

The most important formula for our work is the following. Consider
a function  $\s$ over $\R^\n$. Suppose that $\g$ is admissible. It
thus has a reconstruction wavelet $\r$. An explicit formula has been
given in \ref{recwaveform}.  Then, since
$\Wavetrans_\h=(\Wavetrans_\h\InvWavetrans_\r)\Wavetrans_\g$, the
wavelet transform of $\s$ with respect
 to $\g$  and the one with respect to
$\h$ are related via the so called cross kernel
equation.
\begin{equation}\label{crosskernel}
\Wavetrans_\h\s = \Pi_{\g\to\h}\ast\Wavetrans_\g\s,
\quad \Pi_{\g\to\h} = \Wavetrans_\h\r.
\end{equation}

\section{Some function spaces.}
All formulas that we have given so far have a well defined meaning
if the wavelets are taken in some subset of the class of Schwarz as
we will recall now.

\subsection{The analysis of $\Schwarz_0(\R^\n)$}
Let $\Schwarz(\R^\n)$    denote the class of  Schwarz consisting
of those functions that   together with their derivatives decay faster
 than any polynomial such that the following norms are all finite
for all multi-indices $\alpha$ and $\beta$
$$
\norm{\s}_{\alpha,\beta}=
\sup_{\t\in\R^\n}\abs{\t^\alpha\partial^\beta\s(\t)}<\infty.
$$
They generate a locally convex topology  which makes $\Schwarz(\R^\n)$ a
Fr\'echet space.
 We denoted by $\Schwarz_0(\R^\n)$ the closed set of functions
in $\Schwarz(\R^\n)$ for which all moments vanish
$$
\forall\alpha\in\N^\n, \ \int\d\x\,\s(\x)\, \x^\alpha = 0 \iff
\forall\m>0,\ \fou\s(\k) = \o(\k^\m) \ (\abs\k\to0).
$$
The Schwarz space of functions over the half-plane $\Halfplane^\n$
 shall be denoted by $\Schwarz(\Halfplane)$.
 It consists of those functions $\r$ for which the following norms
 are all finite
$$
\norm{\r}_{\k,\l,\m, \n} =  \sup_{(\b,\a)\in\Halfplane^\n}
\abs{(\a+1/\a)^\k(1+\abs\b)^\l\partial^\m_\b\partial^\n_\a\, \r(\b,\a)}
<\infty.
$$
Note that this means that $\r$ together with all its derivatives decays
 rapidly for large $\b$ and for large or small $\a$. It
 can be shown that
$$
\Wavetrans : \Schwarz_0(\R^\n) \times\Schwarz_0(\R^\n)
\to \Schwarz(\Halfplane^\n),\quad(\g,\s)
\mapsto\Wavetrans_\g\s
$$
is continuous. The same holds for the  wavelet synthesis defined through
$$
\InvWavetrans_\h\r(\x) = \int_{\Halfplane^\n}{\d\b\d\a\over\a}
\, \r(\b,\a) \, {1\over\a^\n}\h\of{\x-\b\over\a},
$$
and we
 have that
$$
\InvWavetrans : \Schwarz_0(\R^\n) \times \Schwarz(\Halfplane^\n)
\to \Schwarz_0(\R),\quad(\h, \r)
\mapsto\InvWavetrans_\h\r
$$
is continuous too. However that in this paper we will not discuss
topologies on the microlocal classes we define. This can be done in
an obvious way and we want to streamline the discussion.

We note here the following important fact. For admissible
$\g\in\Schwarz_0(\R^\n)$, and arbitrary $\h\in\Schwarz_0(\R^\n)$
the crosskernel \ref{crosskernel} is a function in
$\Schwarz(\Halfplane^\n)$. It thus is very well localized.

\subsection{Wavelet analysis of $\Schwarz_0^\prime(\R^\n)$.}

We denote  the space of linear continuous functionals
$\eta: \Schwarz_0(\R^\n)\to\C$ by  $\Schwarz_0^\prime(\R^\n)$.
We consider it together with its natural  weak-$\ast$ topology.
The space $\Schwarz^\prime_0(\R^\n)$ can canonically be identified
with  $\Schwarz^\prime(\R)/\Pol(\R^\n)$, where $\Pol(\R^\n)$ is the
space
 of polynomials in $\n$
variables.
 The wavelet transform of $\eta\in\Schwarz_0^\prime(\R^\n)$
 can now be defined pointwise as
$$
\Wavetrans_\g\eta(\b,\a) = \eta(\conj\g_{\b,\a}).
$$
This is a smooth function (e.g. \cite{HOLSBOOKWAVE}) that satisfies at
\begin{equation}\label{growthconref}
\abs{\Wavetrans_\g\eta(\b,\a)} \leq\const (1+\abs\b)^\m
(\a+\a^{-1})^\m\end{equation}
for some $\m>0$. By duality we have that the mapping (this time for
fixed wavelet $\g\in\Schwarz_0(\R^\n)$)
$$
\Wavetrans_\g :
\Schwarz_0^\prime(\R^\n)\to\Schwarz^\prime(\Halfplane^\n),
\quad
\eta\mapsto\Wavetrans_\g\eta
$$
is continuous. Here $\Schwarz^\prime(\Halfplane^\n)$ is the dual of
$\Schwarz(\Halfplane^\n)$ together with the weak-$\ast$ topology
Vice versa any  let $\r\in\Schwarz^\prime(\Halfplane^\n)$. Then we
may set for $\s\in\Schwarz_0(\R^\n)$
$$
(\InvWavetrans_\h\r)(\s) = \r(\Wavetrans_{\conj\h}\s)
$$
Clearly again ($\h\in\Schwarz_0(\R^\n)$)
$$
\InvWavetrans_\h :
\Schwarz^\prime(\Halfplane^\n)\to\Schwarz_0^\prime(\R^\n),
\quad
\r\mapsto\InvWavetrans_\h\r
$$
is continuous.
In the case of a locally integrable function $\r$ of
at most polynomial growth we have
$$
\InvWavetrans_\h\r(\s) = \int_{\Halfplane^\n}{\d\b\d\a\over\a}\,
 \r(\b,\a) \, \Wavetrans_{\conj\h}\s(\b,\a).
$$

\label{morclodclaksjdhfjksa}\subsection{More general spaces}
Many  other function spaces can be characterized in terms of
wavelet coefficients.
As a rule, the faster the wavelet coefficients decay, the more the
analyzed function is regular.

We come to the details. For this consider a vector space
$\B(\Halfplane^\n)$
 of
locally integrable functions
$$
\Schwarz(\Halfplane^\n)\subset
\B(\Halfplane^\n)\subset\Schwarz^\prime(\Halfplane^\n).
$$
Suppose in addition that $\B(\Halfplane^\n)$ is invariant under
convolutions
with highly localized kernels
$$
\r\in \Schwarz(\Halfplane^\n), \
\s\in \B(\Halfplane^\n) \quad\Rightarrow
\quad
\r\ast\s,\ \s\ast\r \in\B(\Halfplane^\n).
$$
The convolution of two functions over $\Halfplane^\n$ is defined through
\ref{convdefhlfplsne}.
It then makes sense to pull back $\B(\Halfplane^\n)$ to a vector space of
distributions over $\R^\n$. We shall denote this space of distributions
by $\B(\R^\n)$. It is defined through the following theorem

\theorem{difseeriusdf}{For  a distribution
 $\eta\in\Schwarz_0^\prime(\R^\n)$ the following is equivalent.
\begin{itemize}
\item
There is  a wavelet $\g\in\Schwarz_0(\R^\n)$ which is admissible in
that it
satisfies
$$
\int_0^\infty{\d\a\over\a}\abs{\fou\g(\a\k)}^2 \sim 1.
$$
for which we have
$$
\Wavetrans_\g\eta \in \B(\Halfplane^\n).
$$
\item
For all $\h\in\Schwarz_0(\R^\n)$  we have
$$
\Wavetrans_\h\eta\in\B(\Halfplane^\n).
$$
\end{itemize}
}

\Proof
The passage from $\Wavetrans_\g\eta$ to $\Wavetrans_\h\eta$ is
given by the
highly localized cross kernel \ref{crosskernel}. By definition,
this operation leaves invariant
$\B(\Halfplane^\n)$.\qed

Therefore it makes sense to speak of the space $\B(\R^\n)$ associated
with $\B(\Halfplane^\n)$. It is precisely the space of distributions
for which $\Wavetrans_\g\eta\in\B(\Halfplane^\n)$ where $\g$ is
as given in the theorem.

We still shall need an additional technical assumption on the
spaces $\B(\Halfplane^\n)$: their multiplier algebra should  contain
the bounded functions
$$
\m\in\L\infty(\Halfplane^\n), \s\in\B(\Halfplane^\n) \Rightarrow
\m\cdot\s
\in\B(\Halfplane^\n).
$$
This allows us to define the T\"oplitz operators
$$
\InvWavetrans_\h\m\Wavetrans_\g :\quad \B(\R^\n)\to\B(\R^\n).
$$
For the rest of the paper we refer to these spaces
$\B(\Halfplane^\n)$ satisfying all stated properties and
their pulled back counter part  over $\R^\n$,
$\B(\R^\n)$,  as admissible local regularity spaces.

We end this section with a remark concerning topology.
Suppose that in addition $\B(\Halfplane^\n)$ is a Banach space.
Suppose in addition that it is a Banach lattice
$$
\norm{\ \abs\s\ }_{\B(\Halfplane^\n)} = \norm\s_{\B(\Halfplane^\n)}.
$$
Suppose further that for fixed $\Pi\in\Schwarz(\Halfplane^\n)$ we have
$$
\r\mapsto \Pi\ast\r
$$
is continuous.  Then we can define a norm on $\B(\R^\n)$ which
makes it a Banach space by setting
$$
\norm{\s}_{\B(\R^\n)} = \norm{\Wavetrans_\g\s}_{\B(\Halfplane^\n)}.
$$
This is well defined, since for different wavelets satisfying the
hypothesis of theorem \ref{difseeriusdf} we obtain  equivalent
norms.
There is an easy to verify sufficient condition for
$\B(\Halfplane^\n)$ to be stable under convolution with localized
kernels in case that $\B(\Halfplane^\n)$ is a Banach space. It is
enough  to
find $\K$ and $\const$ such that  for all $\s\in\B(\Halfplane^\n)$
we can estimate
$$
\norm{\s(\alpha\cdot+\beta,\alpha\cdot)}_{\B(\Halfplane^\n)}
\leq \const (\alpha+1/\alpha)^\K(1+\abs\b)^\K.
$$
Indeed,  by a simple change of variables, we may write
$$
\norm{\Pi\ast\s}_{\B(\Halfplane^\n)}
\leq \norm{\s}_{\B(\Halfplane^\n)} \, \const\,
\int_{\Halfplane^\n}{\d\beta\d\alpha\over\alpha}
\Pi(\beta,\alpha)(\alpha+1/\alpha)^\K(1+\abs\beta)^\K.
$$
The last integral is a finite constant.
For the sake of simplicity, we shall not give a detailed discussion
of possible topologies on the microlocal classe we are going to
define in the next chapter.

\subsection{Some examples of local regularity spaces}

Many functions spaces of day to day functional analysis can be
characterized
with this easy concept. Most of them are contained in the following
two scales of spaces.
  Let $\phi:\R_+\to\R_+$ and $\kappa:\R^\n\to\R_+$ be two
 tempered weight functions over $\R_+$ and $\R$ respectively. By this we
understand that they satisfy at
$$
\phi(\a\aa)\leq\O(1)(\a+1/\a)^\n\,\phi(\aa),\quad
\kappa(\b+\bb) \leq \O(1) \, (1+\abs\b)^\n\,\kappa(\bb),
$$
for some $\n>0$. Then the following expressions define norms
for functions over the half-space $1\leq\p,\q\leq\infty$
$$
\norm{\r}  = \left\{\int_0^\infty{\d\a\over\a}\phi(\a)
\left(\intinf\d\b\,\kappa(\b)\,\abs{\r(\b,\a)}^\q\right)\right\}^{\p/\q}
$$
or
$$
\norm{\r}  = \left\{\intinf\d\b\,\kappa(\b)
\left(\int_0^\infty{\d\a\over\a}\phi(\a)\,\abs{\r(\b,\a)}^\q\right)
\right\}^{\p/\q}.
$$
The associated Banach spaces are stable under convolution with
highly regular kernels and thus they may be pulled back to $\R^\n$
giving rise to two scales of spaces.
The first scale of spaces contains the Besov spaces, whereas the
second scale contains the $\L\p$-spaces and Sobolev spaces (see e.g.
\cite{MEYERIAII},
\cite{COIFMANMEYER}, \cite{HOLSBOOKWAVE},  \cite{TRIEBELII} for
more details).

For the moment we only note that the space of locally integrable
functions in $\Halfplane^\n$ for which ($\alpha\in\R$)
$$
\abs{\r(\b,\a)}\leq\const\a^\alpha
$$
is stable under convolution with kernels in
$\Schwarz(\Halfplane^\n)$. It can thus be pulled pack to a space of
distributions over $\R^\n$. As is well known by now, these
regularity spaces correspond to the \H\"older---for $\alpha>0$,
$\alpha\not\in\N$---respectively Zygmund
classes---for $\alpha\in\N$. We shall denote this space by
$\Lambda^\alpha(\R^\n)$.

\section{Some more  general microlocal classes.}

In this section we shall be concerned with the problem of
constructing  new  local regularity spaces out of old ones. The idea
is easily explained.
Consider an arbitrary set  $\Omega\subset\Halfplane^\n$. Eventually
we want to consider lines
$\Omega=\{(\b=\lambda\xi,\a=\lambda^\gamma)\}$ in order to define
regularity in the direction $\xi\in\R^\n$. But for the moment we
stay general
since it renders the discussion more easy to understand.
Now fix  in addition to $\Omega$, an admissible,  local regularity
space  $\B(\R^\n)$. It is characterized by the fact that the wavelet
transforms of
 its members are in some vector  space $\B(\Halfplane^\n)$.
 The regularity  classes we want to construct
are roughly speaking distributions whose wavelet coefficients have
on $\Omega$
a growth  behavior governed by the class $\B(\Halfplane^\n)$.

A slightly more general concept is obtained if we take two local
regularity
 spaces $\B_1(\R^\n)$ and $\B_2(\R^\n)$,
 both of the type considered in section \ref{morclodclaksjdhfjksa}.
 We now want to cut the half-space $\Halfplane^\n$ into two parts,
 say $\Omega$
 and its complement $\Omega^{c}$. In some sense---to be made
 precise below---we consider classes of distributions whose wavelet
 coefficients behave inside $\Omega$ like the wavelet coefficients of
 functions in $\B_1(\R^\n)$, whereas in $\Omega^{c}$, these
distributions have regularity governed by $\B_2(\R^\n)$.

A naive approach might be to require that the restriction of
$\Wavetrans_\g\eta$ to $\Omega$ satisfies
$\Wavetrans_\g\eta\in\B_1(\Halfplane^\n)$,  whereas its
restriction to the complement of $\Omega$ should
 correspond (via $\InvWavetrans_\h$)
to a function in $\B_2(\R^\n)$. However
this definition might depend on the wavelets we use  and thus it is
not useful.

To get around this difficulty we first construct a suitable family
of
 neighborhoods for $\Omega$.
With these neighborhoods it
turns out that we can define  vector  spaces that are  independent
of the wavelets $\g$ and $\h$.

\subsection{A non-Euclidian distance}

The first step of the construction is to introduce a non-euclidian
distance function adapted to the geometry of the half-space.  A
suitable choice is given by  $\hbox{dist}$ as defined via
$$
\hbox{dist}((\b,\a), (\bb,\aa))
= \abs{\a/\aa} + \abs{\aa/\a}  + \abs{(\b-\bb)/\aa}
+ \abs{(\bb-\b)/\a}.
$$
This clearly is not a distance  in the usual sense, since
$\hbox{dist}((\b,\a), (\b,\a))=2$. However a kind of multiplicative
triangular inequality holds (see lemma \ref{treisjdhfinghdjska}
below).
Note however that the distance function is  symmetric.
$$
\hbox{dist}((\b,\a), (\bb,\aa))
=\hbox{dist}((\bb,\aa),(\b,\a)).
$$

Clearly the upper half-space carries a natural group structure. It
is given by the following composition law
$$
(\b,\a) (\bb,\aa) = \delta_\a\tau_\b(\bb,\aa) = (\a\bb +\b,\a\aa)
$$
where $\tau_\b$ and $\delta_\a$ stand for the translation and
dilation as left actions of $\R_+$ and $\R^\n$ on $\Halfplane^\n$
via
\begin{equation}\label{tranlsationaction}
\delta_\alpha: (\b,\a) \mapsto (\alpha \b, \alpha \a),\quad
\tau_\beta: (\b,\a) \mapsto
 (\b+\beta,\aa),
\end{equation}
The inverse element  reads
$$
(\b,\a)^{-1} = (-\b/\a,1/\a),
$$
and the neutral element clearly is $(0,1)$.

If we denote by
$\Delta((\b,\a))$ the distance of a point $(\b,\a)$ from the point
$(0,1)$
$$
\Delta((\b,\a)) = \hbox{dist}((\b,\a),(0,1)) = {\a +1/\a}  +
\abs{\b}(1 + 1/a)
$$
then we have the following relation
$$
\hbox{dist}((\b,\a),(\bb,\aa)) = \Delta((\b,\a)^{-1}(\bb,\aa)).
$$
Note also the following identity
$$
\Delta((\b,\a)^{-1}) = \Delta((\b,\a)).
$$
The next lemma shows that a kind of triangular inequality holds

\lemma{treisjdhfinghdjska}{We  have the following triangular inequalities
\begin{eqnarray*}
&\max\{&\Delta((\b,\a))/\Delta((\bb,\aa)),
\Delta((\bb,\aa))/\Delta((\b,\a))\}\\
&&\leq \Delta((\b,\a)(\bb,\aa)) \leq   \Delta((\b,\a))
\Delta((\bb,\aa)).
\end{eqnarray*}
}

\Proof
To prove this inequality note that an elementary direct computation
shows that
\begin{eqnarray*}
\a\aa\Delta((\b,\a)^{-1}(\bb,\aa))  &=&
\a\aa\Delta(((\bb-b)/\a, \aa/\a))\\
&=& \a^2 +\aa^2 +\a\abs{\b-\bb} + \aa \abs{\b-\bb}
\\
&\leq&
\a^2 +\aa^2 +\a\abs{\b} + \a\abs{\bb} + \aa \abs{\b} + \aa\abs{\bb}.
\end{eqnarray*}
On the other hand
$$
\a\aa\Delta((\b,\a))
\Delta((\bb,\aa)) = (1+\a^2 +\abs\b+\a\abs\b)(1+\aa^2
+\abs\bb+\aa\abs\bb).
$$
Therefore the difference between the last  and the previous expression is
majorized by
$$
(1-\a+\a^2)\abs\bb + (1-\aa+\aa^2)\abs\b >0.
$$
The proof of the right most inequality follows now from the identity$$
\Delta((\b,\a)^{-1}) =
\Delta((\b,\a)).
$$

The remaining inequality follows as usual from the previous one, namely$$
\Delta((\b,\a)) \leq  \Delta((\b,\a)(\bb,\aa))\Delta((\bb,\aa)^{-1})
=
\Delta((\b,\a)(\bb,\aa))\Delta((\bb,\aa))
$$
\qed

This immediately  implies the following relation for the distance
function
$$
\hbox{dist}((\b,\a),(\bb,\aa)) = \Delta((\b,\a)^{-1}(\bb,\aa))
\leq \Delta((\b,\a)) \Delta((\bb,\aa))
$$
and
$$
\hbox{dist}((\b,\a),(\bb,\aa))  \geq
\Delta((\b,\a)) /  \Delta((\bb,\aa))
$$
Thus the following triangular inequality holds
$$
\hbox{dist}((\b,\a),(\bbb,\aaa)) \leq
\hbox{dist}((\b,\a),(\bb,\aa)) \hbox{dist}((\bb,\aa),(\bbb,\aaa))
$$
and
$$
\hbox{dist}((\b,\a),(\bbb,\aaa)) \geq
\hbox{dist}((\b,\a),(\bb,\aa)) /\hbox{dist}((\bb,\aa),(\bbb,\aaa)).
$$

\subsection{A family of neighborhoods.}
Let us introduce the closed non-euclidian balls
$$
\U((\b,\a), \r) = \{(\bb,\aa)\in\Halfplane \suchthat
\hbox{dist}((\b,\a), (\bb,\aa))\leq \r\}
$$
Note that they all are obtained by dilations and translations of
the balls around the point $(0,1)$. More precisely, since the
distance function satisfies at
$$
\hbox{dist}((\gamma\b+\beta, \gamma\a), (\gamma\bb + \beta,
\gamma\aa))= \hbox{dist} ((\b,\a), (\bb,\aa)), \quad
\gamma>0, \beta\in\R^\n,
$$
we have
$$
\U((\b,\a), \r) = \tau_\b\delta_\a \U((0,1), \r).
$$
An equivalent system of neighborhoods $\U^\prime$ is obtained by
translating and dilating the family of balls defined via the
following inequalities
$$
(\a-[1+\r+1/(1+\r)]/2)^2 + \abs\b^2 \leq
(1+\r+1/(1+\r))^2.
$$
They are euclidian balls with the "south-pole" at the point
$(\b=0,\a=1/(1+\r))$
and the \lq\lq north-pole\rq\rq\ at the point $(\b=0,\a=(1+\r))$.
The equivalence being expressed by the fact that for some constants
$\const>1$   we have
$$
\U^\prime((0,1), \r/\const) \subset \U((0,1), \r)\subset
\U^\prime((0,1), \const\r).
$$
We may leave the elementary calculations to the reader.

We now are interested in when it makes sense to speak of a certain
regularity in one set and  an other regularity in an other set of
the half-space.
Consider therefore two arbitrary subsets $\Omega$,
$\Sigma\subset\Halfplane^\n$.  We now say that $\Omega$ and $\Sigma$
are well separated  if the following holds.
For  $(\b,\a)\in\Omega$  consider a non-euclidian ball $\U((\b,\a),\r)$
with center $(\b,\a)$ and radius $\r$. Choose $\r$ small enough so
that $\U$ does not meet $\Sigma$. Well separated means for us that
for
some $\epsilon>0$   we may choose $\r$ such that the following
estimate holds true for small $\a$
$$
\r > \a^{-\epsilon}.
$$
In other terms we define  more formally and slightly more general

\definition{}{We say that two sets $\Omega$ and $\Sigma$ are well
separated if  for some  $\epsilon>0$   we have for
\begin{equation}\label{wellsepdrer}
(\b,\a)\in\Omega\Rightarrow
\hbox{dist}((\b,\a), \Sigma) > \Delta((\b,\a))^\epsilon.
\end{equation}
}

Here the distance between a point and a set
$\Omega\subset\Halfplane^\n$ is defined
as usual
$$
\hbox{dist} ((\b,\a), \Omega) = \inf_{(\bb,\aa)\in\Omega}
\hbox{dist} ((\b,\a), (\bb,\aa)).
$$

Note that although the non-euclidian distance diverges at small
scale,the euclidian distance might tend to $0$ as $\a\to0$.
As an example, which is somehow typical,  consider in
$\Halfplane^2$ the sets
$$
\Omega=\{\a>\abs\b^\alpha\}\cap\{\a<1/2\},\quad
\Sigma=\{\a<\abs\b^\beta\}\cap\{\a<1/2\}.
$$
They are well separated iff $\alpha>\beta$. However their euclidian
distance tends always to $0$ as $\a\to0$.

We have still an other useful characterization of well separatedness of
two sets $\Omega$ and $\Sigma$. For this consider the sets
$$
\delta_{1/\a}\tau_{-\b}\Omega,\quad(\b,\a)\in\Sigma.
$$
Now both sets are well separated iff  each  of these sets is
contained in the complement of
 a non-Euclidian ball $\U((0,1),\r(\b,\a))$, with $\r(\b,\a) >
\Delta((\b,\a))^{\epsilon}$.
\begin{equation}\label{wellsepdacharz}
\delta_{1/\a}\tau_{-\b}\Omega\subset \U((0,1),
\Delta((\b,\a))^{\epsilon})^{c}
\end{equation}

\subsection{More about well separated sets.}

Since the distance function is continuous in the euclidian
topology, it is clear that the distance of a point and a set and its
euclidian closure are the same.
Therefore a set is well separated form an other if and only if its
closures
are.

The notion of well separated is inherited by subsets. If
$\Omega\subset\Sigma$ and $\Sigma$ is well separated from
$\Xi$, then  $\Omega$ is well separated from $\Xi$ too.

The notion of well separated is symmetric and we have

\lemma{symnetriocidusio}{If $\Omega$ is well separated from
$\Sigma$  then $\Sigma$ is well separated from $\Omega$.
}

\Proof
By hypothesis we have that $(\b,\a)\in\Omega$ and $(\bb,\aa)\in\Sigma$
implies that
\begin{equation}\label{wellresieufjsa}
\hbox{dist}((\b,\a),(\bb,\aa)) >  \Delta((\b,\a))^{\epsilon}
\end{equation}
We claim that this implies that
$$
\dist((\b,\a),(\bb,\aa)) >   \Delta((\bb,\aa))^{\epsilon^\prime}
$$
for $\epsilon^\prime=\epsilon/(1+\epsilon)$. For suppose that on
the contrary for some points we have
$$
\dist((\b,\a),(\bb,\aa))  = \Delta((\b,\a)^{-1} (\bb,\aa))
\leq \Delta((\bb,\aa))^{\epsilon^\prime}
$$
This implies, via the second triangular inequality,  in particular that$$
\Delta((\bb,\aa))/ \Delta((\b,\a)) \leq
\Delta((\bb,\aa))^{\epsilon^\prime}
$$
and therefore by the choice of $\epsilon^\prime$ that
$\Delta((\bb,\aa))^{\epsilon^\prime}\leq
\Delta((\b,\a))^{\epsilon}$. This implies
$$
\hbox{dist}((\b,\a),(\bb,\aa))
\leq \Delta((\b,\a))^{\epsilon},
$$
which  is in  contradiction with \ref{wellresieufjsa}.
\qed

For every $\epsilon>0$ let us introduce the following non-Euclidian
neighborhoods of a set $\Omega\subset\Halfplane^\n$
$$
\Gamma_\epsilon(\Omega) = \bigcup_{(\b,\a)\in\Omega}
\U\big((\b,\a), \ \Delta((\b,\a))^{\epsilon}\big).
$$
The system of such neighborhoods constitutes a fundamental family
of neighborhoods in the following sense. We have that
$\Gamma_\epsilon(\Omega)^{c}$ is well separated
from $\Omega$. In addition, if $\Sigma$
is well separated from $\Omega$, then for some $\epsilon>0$   we
have that
$$
\Sigma\cap \Gamma_\epsilon(\Omega) = \phi.
$$

Thanks to the triangular inequalities we have the following
associativity for the $\epsilon$ neighborhoods.

\lemma{weciteitnow}{For any set $\Omega\subset\Halfplane^\n$ the
following holds true. For $\epsilon_1>0$, $\epsilon_2>0$ and
$\epsilon_3\geq
\epsilon_1+\epsilon_2(1+\epsilon_1)$
we have
$$
\Gamma_{\epsilon_2}(\Gamma_{\epsilon_1}(\Omega)) \subset
\Gamma_{\epsilon_3}(\Omega).
$$
On the other hand for $\epsilon_2$ such that
$\epsilon_2/(1-\epsilon_2) <\epsilon_1$, there is an $\epsilon_3>0$
such that
$$
\Gamma_{\epsilon_2}(\Gamma_{\epsilon_1}(\Omega)^{c})^{c} \supset
\Gamma_{\epsilon_3}(\Omega).
$$
More precisely  it is enough that the following relation holds
$$
\epsilon_1 - \epsilon_3 - {\epsilon_2(1+\epsilon_3)\over(1-\epsilon_2)}
\geq0.
$$
}

\Proof
We show the first part.
If $(\bbb,\aaa)\in\Gamma_{\epsilon_2}(\Gamma_{\epsilon_1}(\Omega))$ then
for some points $(\bb,\aa)\in\Halfplane^\n$ and $(\b,\a)\in\Omega$
we have
$$
\hbox{dist}((\bbb,\aaa),(\bb,\aa)) \leq
\Delta((\bb,\aa))^{\epsilon_2},\quad
\hbox{dist}((\bb,\aa),(\b,\a)) \leq \Delta((\b,\a))^{\epsilon_1}.
$$
Therefore by the triangular relation we have
$$
\hbox{dist}((\bbb,\aaa),(\b,\a)) \leq
\Delta((\bb,\aa))^{\epsilon_2}\Delta((\b,\a))^{\epsilon_1}.
$$
Now as before by the reverse triangular inequality we have
$$
\Delta((\bb,\aa))/\Delta((\b,\a))\leq \Delta((\b,\a))^{\epsilon_1}
$$
 and therefore  finally as claimed
$$
\hbox{dist}((\bbb,\aaa),(\b,\a))
\leq \a^{-\epsilon_1-\epsilon_2(1+\epsilon_1)},
$$

The   second statement may be rephrased as follows:   if for all
$(\b,\a)\in\Omega$ we have
\begin{equation}\label{fhdufydusiaudyf}
\hbox{dist}((\b,\a),(\bbb,\aaa)) \leq \Delta((\b,\a))^{\epsilon_3}
\end{equation}
then
$(\bbb,\aaa)\not\in\Gamma_{\epsilon_2}(\Gamma_{\epsilon_1}(\Omega)^{c})$.
Suppose that the contrary is true. Then for some
$(\bb,\aa)\in\Halfplane^\n$ satisfying at
$$
\forall(\b,\a)\in\Omega \quad
\hbox{dist}((\b,\a),(\bb,\aa))> \Delta((\b,\a))^{\epsilon_1},
$$
we have
\begin{equation}\label{scofndhjfy}
\hbox{dist}((\bbb,\aaa),(\bb,\aa))\leq \Delta((\b,\a))^{\epsilon_2}.
\end{equation}
Therefore we have by the triangular inequality
\begin{eqnarray*}
\Delta((\b,\a))^{\epsilon_1}&<&\hbox{dist}((\b,\a),(\bb,\aa))\\
&\leq&\hbox{dist}((\b,\a),(\bbb,\aaa))
\hbox{dist}((\bbb,\aaa),(\bb,\aa))\\
&\leq &
 \Delta((\b,\a))^{\epsilon_3}\Delta((\bb,\aa))^{\epsilon_2}.
\end{eqnarray*}
Now, by the reverse triangular inequality, \ref{fhdufydusiaudyf} implies
$$
\Delta((\bbb,\aaa))/\Delta((\b,\a))\leq\Delta((\b,\a))^{\epsilon_3}.
$$
Now again by the reverse triangular inequality \ref{scofndhjfy}
implies
$$
\Delta((\bb,\aa))/
\Delta((\bbb,\aaa))\leq\Delta((\bb,\aa))^{\epsilon_2}.
$$
It follows that
$$
\Delta((\bb,\aa))\geq \Delta((\b,\a))^{(1+\epsilon_3)/(1-\epsilon_2)}.
$$
Therefore
$$
1 < \Delta((\b,\a))^{\epsilon_3 -\epsilon_1 + \epsilon_2
(1+\epsilon_3)/(1-\epsilon_2)},
$$
which is impossible by the choice of $\epsilon_3$ and  since
$\Delta((\b,\a))\geq2$.
\qed

This immediately implies the following corollary that we shall use
in the next section.

\lemma{wedsindtermed}{Let $\Sigma\supset\Omega$ be such that
$\Sigma^{c}$ is well separated from $\Omega$. Then there is a set
$\Xi$, $\Sigma\supset\Xi\supset\Omega$
such that $\Xi$ is well separated from $\Sigma^{c}$ and $\Omega$ is
well separated from $\Xi^{c}$.}

\Proof
Some $\Gamma_\epsilon(\Omega)$ with $\epsilon$ small enough will do.
\qed

\subsection{Cutting the half-space.}
Let us come back to our original goal of dividing the half-space
into two set of different regularity. As we said already, it is not
possible to speak
of regularity $\B(\Halfplane^\n)$ inside a given set
$\Omega\subset\Halfplane^\n$, since notion is not independent under
highly regular Calder\'on Zygmund operators, or to put it simpler,
it might depend on the given wavelet we use for the definition.

However if we require regularity $\B(\Halfplane^\n)$ in a region
that that is slightly larger than $\Omega$, it then follows that the
same regularity holds true in $\Omega$ for any wavelet.

 We denote by abuse of notation $\Sigma$, (respectively $\Omega$)
the operator that restricts functions over
$\Halfplane^\n$ to the set $\Sigma$ (respectively $\Omega$). That
is we have
$$
\Sigma  : \r\mapsto \chi_{\Sigma}\r,
$$
where $\chi_{\Sigma}$ is the characteristic function of
$\Sigma$.

\theorem{fhduwitheifdhs}{ Consider two  sets $\Sigma$ and $\Omega$
and suppose that $\Sigma\supset\Omega$ in such a way that
$\Sigma^{c}$ and $\Omega$ are well separated.  Let $\g$, $\gg$,
$\h$, $\hh\in\Schwarz_0(\R^\n)$ satisfy  for some $\const>1$
($\s=\g$, $\gg$, $\h$, $\hh$),
$$
\const^{-1} < \int_0^\infty{\d\a\over\a}\abs{\fou\s(\a\k)}^2<\const.
$$
Suppose that $\eta\in\Schwarz^\prime(\R^\n)$ satisfies at
$$
\InvWavetrans_\h\Sigma\Wavetrans_\g\eta \in\B(\R^\n)
$$
Then
$$
\InvWavetrans_{\hh}\Omega\Wavetrans_\gg\eta  \in\B(\R^\n).
$$
}

Therefore  it makes sense to separate $\Halfplane^\n$ into regions
of different regularity provided, the regions are well
 separated.

The proof is based on the following lemma. It estimates the
influence under
 convolution operators over the half-plane of a nasty function
inside some
 region $\Xi$ on $\Xi^\prime$ when both are well separated.

\lemma{}{Suppose $\Xi^\prime$ and $\Xi$ are well separated.
Let $\r$ be a locally integrable function over $\Halfplane^\n$
 that is equal to $0$ except on $\Xi^\prime$, where it satisfies
 for some $\M>0$ and some $\const>0$
$$
(\b,\a)\in\Xi^\prime\quad \Rightarrow\quad
\abs{\r(\b,\a)}\leq\const\,\Delta((\b,\a))^\M.
$$
Then for $\Pi\in\Schwarz(\Halfplane^\n)$ we have that
$\tilde\r=\Pi\ast\r$ satisfies at
$$
(\b,\a)\in\Xi\quad \Rightarrow\quad
 \abs{\tilde\r(\b,\a)}\leq
 \const_\k\Delta((\b,\a))^{-\k}
$$
for all $\k\in\N$.}

\Proof
 We have  to estimate the localization of $\Pi\ast\r(\b,\a)$ for
 $(\b,\a) \in\Xi$.
By definition $\Pi\ast\r(\b,\a)$ equals
$$
\int_{\Xi^\prime}{\d\bb\d\aa\over\aa}
{1\over\aa^\n}\Pi\of{{\b-\bb\over\aa},{\a\over\aa}}\r(\bb,\aa)
$$
By dilation and translation we also may write using the action
(\ref{tranlsationaction}) of dilation and translation on
$\Halfplane^\n$
$$
\int_{\delta_{1/\a}\tau_{-\b}\Xi^\prime}{\d\bb\d\aa\over\aa}
\Pi(\bb,\aa) \, \r(\aa\bb+\b, \a\aa).
$$
Now by hypothesis on $\r$ we may write with some $\K$  and some
$\const>0$
$$
\abs{\r(\aa\bb+\b, \a\aa)}
\leq \const \, \Delta((\b,\a))^\K\Delta((\bb,\aa))^\K.
$$
Plugging this estimation into the previous expression we have to estimate
for $(\b,\a)\in\Xi$
\begin{equation}\label{toeidufghdsja}
\Delta((\b,\a))^\K
\int_{\delta_{1/\a}\tau_{-\b}\Xi^\prime}{\d\bb\d\aa\over\aa}
\abs{\Pi^\prime(\bb,\aa)}
\end{equation}
with $\Pi^\prime(\bb,\aa) = \Delta((\bb,\aa))^\K\Pi(\bb,\aa)$.
Together with $\Pi$  we have that $\Pi^\prime$ is highly localized.
For $\lambda\geq0$ let us look at the following integral running
over the complement of a non-euclidian ball centered at $(0,1)$
$$
\F(\lambda) = \int_{\Delta((\bb,\aa))>\lambda}
{\d\bb\d\aa\over\aa}\,\abs{\Pi^\prime(\bb,\aa)}.
$$
Thanks to the high localization of $\Pi^\prime$, this function is
faster decaying than any power of $\lambda$ as $\lambda\to0$.  Now,
since $\Xi$ and $\Xi^\prime$ are well separated we may use
characterization \ref{wellsepdacharz}
to conclude that the integral in \ref{toeidufghdsja} is estimated
by $\F(\Delta((\b,\a))^\epsilon)$ for some $\epsilon>0$. But this
function is again rapidly decaying as $\Delta((\b,\a))\to\infty$ and
the proof is finished.
\qed

Note that the lemma we went to prove may be rephrased as follows:
for all $\Pi\in\Schwarz(\Halfplane^\n)$ we have that
$$
\s\mapsto \Xi^\prime\,(\Pi\ast(\Xi\s)))
$$
is infinitely  smoothing in the sense that it maps functions of
polynomial growth into rapidly decaying functions over the
half-space.

\Proof (of theorem \ref{fhduwitheifdhs})
These previous considerations  imply   the following: if we have
that $\Sigma\Wavetrans_\g\eta
\in\B(\Halfplane^\n)$, with $\g$ admissible, then for all
$\gg\in\Schwarz_0(\R)$ we have that
$\Omega\Wavetrans_\gg\eta\in\B(\Halfplane^\n)$. To show this note
the passage from $\Wavetrans_\g\eta$ to $\Wavetrans_\gg\eta$ is done
by convolving with
a highly localized kernel $\Pi$. Now we may write
$$
\Omega\Wavetrans_\gg\eta = \Omega(\Pi\ast(\Sigma\Wavetrans_\g\eta))
+ \Omega(\Pi\ast(\Sigma^{c}\Wavetrans_\g\eta)).
$$
Since by hypothesis $\B(\Halfplane^\n)$ is invariant under
multiplications with bounded functions and convolutions with $\Pi$
the first term is again in $\B(\Halfplane^\n)$, whereas the second
term is arbitrary smooth.

A slightly more complicated situation occurs in our theorem, since
we can not conclude from $\InvWavetrans_\h\r\in\B(\R^\n)$ that
$\r\in\B(\Halfplane^\n)$ since the wavelet synthesis is not
injective.

Now we can find a set $\Xi$ between $\Sigma$ and $\Omega$
$$
\Omega\subset\Xi\subset\Sigma,
$$
such that $\Xi$ is well separated from the complement of $\Sigma$
and$\Omega$ is well separated from the complement of $\Xi$. This
follows from lemma \ref{wedsindtermed}. We may conclude that
$$
\Xi\Wavetrans_\g\InvWavetrans_\h\Sigma\Wavetrans_\g\eta
= \Xi (\Pi_1 \ast(\Sigma\Wavetrans_\g\eta) )\in\B(\Halfplane^\n).
$$
Where $\Pi_1=\Wavetrans_\f\h$ for any admissible
$\f\in\Schwarz_0(\R^\n)$.
In particular we may choose $\f$ to be a reconstruction wavelet for
$\g$ and thus it follows  that
$\Pi_1 \ast \Wavetrans_\g\eta = \Wavetrans_\g\eta$.
Now writing (as characteristic functions!) $\Sigma=1-\Sigma^{c}$
the last expression
 equals
$$
\Xi (\Pi_1 \ast\Wavetrans_\g\eta)
 - \Xi (\Pi_1 \ast(\Sigma^{c}\Wavetrans_\g\eta)).
$$
The  set $\Xi$ is well separated from $\Sigma^{c}$ and thus
the second term has rapid decay as $\Delta((\b,\a))$ gets large.
Let us call this function $\u$. Then since $\u$ is well localized we
have
$$
\r\ast\u\in\Schwarz(\Halfplane^\n)
$$
for all $\r\in\Schwarz(\Halfplane^\n)$.
We therefore obtain, up to a function of rappid decay
$$
\Xi \Wavetrans_\g\eta  \in\B(\Halfplane^\n).
$$
 Now $\Wavetrans_\gg\eta=\Pi\ast\Wavetrans_\g\eta$ for some
$\Pi\in\Schwarz(\Halfplane^\n)$. Therefore, since $\Omega\subset\Xi$
is well separated from the complement of $\Xi$ we have as at the
beginning of the proof
that  $\Omega\Wavetrans_\gg\eta\in\B(\Halfplane^\n)$ up to the well
localized function $\Omega\u$.
But then clearly $\InvWavetrans_\hh\Omega\Wavetrans_\gg\eta\in\B(\R^\n)$.
\qed

Let $\Sigma\supset\Omega$ be open and let again $\Omega$ be well
separated from the complement of $\Sigma$ as before. Consider the
two T\"oplitz operators
$$
\T_\Sigma = \InvWavetrans_\h\Sigma\Wavetrans_\g,\quad
\T_\Omega = \InvWavetrans_\h\Omega\Wavetrans_\g.
$$
We then have proved the following
\corollary{}{We have that
$$
[\T_\Sigma,\T_\Omega] = \T_\Sigma \T_\Omega - \T_\Omega\T_\Sigma
$$
is infinitely smoothing in the sense that it maps the tempered
distributions in
$\Schwarz_0^\prime(\R^\n)$ into smooth function in $\Schwarz_0(\R^\n)$.}

\subsection{Some microlocal classes.}

The theorems of the previous section may be used to define some
very general micro-local classes.
Suppose we are given two regularity spaces $\A(\R^\n)$ and
$\B(\R^\n)$ and  Suppose in addition that $\B(\R^\n)\subset
\A(\R^\n)$. Consider a set  $\Omega\subset\Halfplane^\n$. Since we
are only interested in local properties we may suppose that $\Omega$
is bounded in the euclidian norm. In order to avoid technicalities
we suppose $\Omega$  is closed.

  The first type of local regularity classes corresponds to the
idea the globally, a distribution has a regularity discribed by
$\A(\Halfplane^\n)$ whereas locally, in $\Omega$ we have some higher
regularity of type $\B(\Halfplane^\n)$.

A dual idea would be to have the wavelet coefficients concentrated
on the subset $\Omega$. That is that outside of $\Omega$, the
wavelet coefficients are small, hence correspond to the higher
regularity $\B(\Halfplane^\n)$, whereas inside $\Omega$, the
coefficients are in $\A(\Halfplane^\n)$.

We now want to make these statements more precise.
Consider first the case of higher local regularity.
 Suppose that there is  a sequence of closed sets $\{\Omega_\k\}$,
$\k=1,2,\dots$ with
$$
\Omega_1\subset \dots \Omega_\k\subset\Omega_{\k+1}\dots\subset\Omega.
$$
We suppose that $\Omega_\k$ converges to $\Omega$ in the sense that
$$
\Omega = \bigcup_{\k}\Omega_\k.
$$
Suppose that   $\Omega_\k^{c}$ and $\Omega_{\k+1}$ are well
separated for each $\k$. Then clearly $\Omega_\k^{c}$ and
$\Omega_\l$
are well separated for $\k< \l$.  We then say that
$\eta\in\Schwarz_0^\prime(\R^\n)$ belongs to
the microlocal class $\Omega_{\A,\B}$ iff for some admissible $\g$
and all $\k$ we have
$$
\eta\in\A(\R^\n),\quad\hbox{and}\quad
\InvWavetrans_\g\Omega_\k\Wavetrans_\g\eta\in\B(\R^\n).
$$
By the results of the previous theorem
 it is clear that the definition does not
depend on the specific wavelets nor on the family of approximating
sets $\Omega_\k$.  Indeed, by lemma \ref{weciteitnow}   we may take
the family $\Gamma_{1/\k}(\Omega)$ as universal family of
approximating sets.

Note however that for arbitrary $\Omega$, the previous class might
coincide with $\A(\R^\n)$. Indeed, in order to have an approximating
sequence from the interior, of mutually well separated sets the
\lq\lq smoother\rq\rq\ region can not be arbitrary thin. It must
contain at least some non-euclidian neighborhood of some set.

Frequently one takes $\A=\Schwarz^\prime_0(\R^\n)$ in which case
one is only interested in the behavior of the wavelet coefficients
around $\Omega$. In the next section we shall
use this kind of classes to define directional regularity in
distributions.

Consider now the  dual approach, where we want to formalize the
idea of wavelet coefficients concentrated on $\Omega$.  Suppose now
that a sequence of open sets $\{\Omega_\k\}$, $\k=1,2,\dots$ with
$$
\Omega_1\supset \dots \Omega_\k\supset\Omega_{\k+1}\dots\supset\Omega.
$$
 converges to $\Omega$ in the sense that
$$
\Omega = \bigcap_{\k}\Omega_\k.
$$
Again we require that   $\Omega_\k^{c}$ and $\Omega_{\k+1}$ are
well separated for each $\k$.  We then say that
$\eta\in\Schwarz_0^\prime(\R^\n)$ belongs to
the microlocal class $\Omega^{\A,\B}$ iff for some admissible $\g$
and all $\k$ we have
$$
\eta\in\A(\R^\n),\quad\hbox{and}\quad
\InvWavetrans_\g\Omega_\k^{c}\Wavetrans_\g\eta\in\B(\R^\n).
$$
Note that in the case where $\B(\R^\n)=\Schwarz_0(\R^\n)$, this
corresponds to
the idea having the wavelet coefficients concentrated on the set
$\Omega$, where they satisfy the less restrictive regularity
estimate given by $\A(\Halfplane^\n)$.

Note again that the region that contains the wavelet coefficients
corresponding to the smoother behavior can not be arbitrary thin.
However the set $\Omega$
 on which the wavelet coefficients are concentrated is arbitrary.

\label{sec: singdirect@wavfst}\section{Some directional microlocal
classes.}

We now propose to look at more specific examples of regularity
classes. In particular to those we mentioned in the beginning of the
paper, that is to classes related to the notion of singular or
regular directions in distributions.

Particular useful examples arise when we consider parabolic regions
or lines
in wavelet space. As measure of regularity it is useful to consider
the H\"older-Zygmund scale $\Lambda^\alpha$ of spaces defined in
wavelet space via
$$
\norm{\s}_\alpha =
\sup_{(\b,\a)\in\Halfplane^\n}\abs{\a^{-\alpha}\s(\b,\a)}
$$
Now fix a  vector $\xi\in\R^\n$, $\abs\xi>0$ and consider the set
$$
\Xi = \Xi(\xi,\gamma) =\{(\b = \lambda\xi, \a =
\lambda^\gamma)\suchthat 1/2>\lambda>0\}
$$ for some $\gamma>1$. We now  say that
$\eta\in\Lambda^\alpha(\R^\n)$ is locally of type
$(\alpha,\xi,\gamma)$ if it belongs to the
microlocal class $\Omega^{\A,\B}$ with
$\Omega=\Gamma_\epsilon(\Xi(\xi,\gamma))$ for some $\epsilon>0$,
$\A=\Schwarz(\R^\n)$ and $\B=\Lambda^\alpha(\R^\n)$.
Explicitly, this means---let us recall it once more---that
the wavelet transform of $\s$ satisfies for some $\epsilon>0$ at
$$
(\b,\a) \in\Gamma_\epsilon(\Xi(\xi,\gamma))\Rightarrow
\abs{\Wavetrans_\g\eta(\b,\a)}\leq \const \a^\alpha.
$$
and
$$
(\b,\a) \not\in\Gamma_\epsilon(\Xi(\xi,\gamma))\Rightarrow
\abs{\Wavetrans_\g\eta(\b,\a)}\leq \const  (\a+1/\a)^{\k}(1+\abs\b)^{\k}.
$$
This corresponds to looking at behavior of the wavelet coefficients
under the following non-homogeneous dilations
$$
\Wavetrans_\g\eta(\lambda\xi,\lambda^\gamma),\quad
\lambda>0.
$$
Here $\const$ depends on $\epsilon$ only.
Actually we may choose $\xi$ such that $\abs\xi=1$. Indeed suppose
$\xi^\prime=\beta\xi$ with $\beta>0$ and denote by $\Xi^\prime$ the
corresponding line. Then if $(\b,\a)\in\Xi$ it follows that
$(\b,\beta^{-\alpha}\a)\in\Xi^\prime$. Therefore the non-euclidian
distance between the two points is uniformly bounded by
$\beta^\alpha+\beta^{-\alpha}$. Therefore, they define the same
micro-local classes.

Let us explain in which sense these classes are linked to singular
and regular directions. For this replace for the moment the wavelet
at position $\b$ and scale $\a$ by the characteristic function of
the euclidian ball centered at $\b$ and of radius $\a$. As $(\b,\a)$
tends to $(0,0)$, while always
$(\b,\a)\in(\Omega(\xi,\gamma))$, the support of the wavelets is
contained in
a cusp-like region, around the line in direction $\xi$. This shows,
thatthe micro-local class  $(\alpha,\xi,\gamma)$ quantifies the
regularity of
$\s$ in direction $\xi$.

\subsection{Some elliptic regularity}

We now want to apply the classes introduced above to a problem of
elliptic regularity. For the sake of simplicity, we only discuss
the Laplace equation and leave more general elliptic operators for
subsequent papers.
We say an open  domain $\Omega\subset\R^\n$ satisfies the cusp
condition of
degree $\delta>0$ at  $\x\in\partial\Omega$ in direction $\xi$,
$\xi\in\R^\n-\{0\}$ if there is  some $\const>0$ such that
$$
\{\y\suchthat \abs{\y-\x}\leq \const \abs{(\y-x\mid\xi)}^\delta,
 \abs{\y-\x}\leq \const \} \subset \Omega
$$

\theorem{}{Let $\Omega$ satisfy a cusp condition at $0$ of type
$\delta$ in direction $\xi$.
Suppose that $\eta$ is a tempered distribution that satisfies
inside $\Omega$ at
$$
\Delta\eta = \f
$$
for some   distribution $\f$ supported by $\conj\Omega$.
Now if $\f$ is of type $(\xi,\gamma,\alpha)$, with $\gamma>\delta$,
then  it follows that
$\eta$ is of type $(\xi,\gamma,\alpha+2)$.}

This theorem is a special case of a more general theorem. Let
$\B(\R^\n)$ be a local regularity space of the type we have
considered before with
$\B(\Halfplane^\n)$ the associated space of functions over
$\Halfplane^\n$.
It is plain to see that together with $\B(\Halfplane^\n)$, the
space of functions which consists the functions $\a^\gamma\r(\b,\a)$
with $\r\in\B(\Halfplane^\n)$ is again an admissible regularity
space.
It shall be denoted by $(\a^\gamma\B)(\Halfplane^\n)$ respectively
$(\a^\gamma\B)(\R^\n)$,

For a set $\Omega\subset\R^\n$ we consider the set
$$
\bigcup_{\b\in\Omega}\K(\b),
$$
where $\K(\b)\subset\Halfplane^\n$ is the cone of opening angle $1$
with top in $\b$
$$
\K(\b) = \{(\beta,\alpha)\suchthat \abs{\beta-\b}\leq \alpha\}.
$$
We call this set the influence region of $\Omega$ in the upper
half-space.

The general theorem can now be stated as follows.

\theorem{}{Let $\Omega\subset\R^\n$ be open. Suppose that
$\Xi\subset\Halfplane^\n$ is well separated from the influence
region of $\Omega$. Suppose that $\Xi^\prime\subset\Xi$ such that
$\Xi^\prime$ and $\Xi^{c}$ are well separated. Suppose that $\eta$
is a tempered distribution that satisfies inside $\Omega$ at
$$
\Delta\eta = \f
$$
for some distribution $\f\in\Schwarz(\R^\n)$ supported by
$\conj\Omega$. If now $\f$ satisfies at
$$
\InvWavetrans_\g\Xi\Wavetrans_\g\f \in \B(\R^\n),
$$
with some   admissible wavelet $\g\in\Schwarz_0(\R^\n)$,
then it follows that
$$
\InvWavetrans_\g\Xi^\prime\Wavetrans_\g\eta \in (\a^2\B)(\R^\n),
$$
}

\Proof
We may suppose that $\g\in\Schwarz_0(\R^\n)$ be spherically
symmetric, $\g=\g(\abs\x)$.
Then with $\h=\Delta\g$ we may write
$$
\Wavetrans_\g\Delta\eta = -\a^{-2}\Wavetrans_\h\eta.
$$
Now $\g$ and $\h$ are both admissible in the sense that
$$
\int_0^\infty{\d\a\over\a}\abs{\fou\g(\a\k)}^2 \sim
\int_0^\infty{\d\a\over\a}\abs{\fou\h(\a\k)}^2 \sim 1
$$
 From this it follows immediately that
$\s\in\Schwarz^\prime(\R^\n)$ and $\Delta\s$ satisfying at
$$
\InvWavetrans_\g\Xi\Wavetrans_\g\Delta\s \in \B(\R^\n),
$$
implies that $\s$ satisfies at
$$
\InvWavetrans_\g\Xi^\prime\Wavetrans_\g\s \in (\a^2\B)(\R^\n).
$$
The theorem is therefore proved if
we can show that any  distribution $\f^\prime$ that coincides with
 $\f$ inside $\Omega$ satisfies again at
$$
\InvWavetrans_\g\Xi\Wavetrans_\g\Delta\f \in \B(\R^\n).
$$
But this is follows from the next lemma, which justifies the name
influence region for $\Omega$.

\lemma{}{Let $\Xi\subset\Halfplane^\n$ be well separated from the
influence region of $\Omega$. Then for all
$\rho\in\Schwarz^\prime(\R^\n)$ with support in $\Omega$ we have
that
$$
\InvWavetrans_\g\Xi\Wavetrans_\g\rho\in\Schwarz_0(\R^\n).
$$
}

\Proof
By hypothesis
there is an $\epsilon>0$ such that for
$(\b,\a)\in\Gamma_\epsilon(\Xi)$ there is a Euclidian ball of radius

$>\Delta((\b,\a))^{\epsilon}$ around $\b$ that is contained in the
complement of the influence region of $\Omega$. Denote by
$\phi$ a $\CC^\infty$ function which is identically $1$ on the
complement of the unit ball of radius $1$ and which is supported on
the complement of a slightlyslightlysmaller ball. Then denote by
$\phi_{\b,\a}$ the family of translates and dilates of $\phi$.
Therefore if $\rho$ is supported by $\Omega$, it follows that
$$
(\b,\a)\in\Xi\Rightarrow
\Wavetrans_\h\rho(\b,\a) =
\Wavetrans_\g(\phi_{\b,\Delta((\b,\a))^{\epsilon}} \rho)(\b,\a)
=\scalarproduct{\g_{\b,\a}\phi_{\b,\Delta((\b,\a))^{\epsilon}}}{\rho}
$$
But now for every $\epsilon>0$ we have that
$\g_{\b,\a}\phi_{\b,\Delta((\b,\a))^{\epsilon}}$ tends to $0$ as
$\Delta((\b,\a))\to\infty$ in $\Schwarz(\R^\n)$ in such a way that
for all semi-norms in $\Schwarz_0(\R^\n)$
and all $\M>0$ we have
$$
\norm{\g_{\b,\a}\phi_{\b,\Delta((\b,\a))^{\epsilon}}}_{\l,\m}\leq
\const_{\l,\m,\M}\Delta((\b,\a))^{-\M}
$$
This proofs the lemma
\qed
The theorem is proved.
\qed


\begin{thebibliography}{Mey90}

\bibitem[CM90]{COIFMANMEYER}
R.~Coifman and Y.~Meyer.
\newblock {\em Ondelettes et Op\`erateurs III}.
\newblock Hermann, 1990.

\bibitem[Dau92]{DAUBECHIESI}
I.~Daubechies.
\newblock {\em Ten Lectures on Wavelets}.
\newblock SIAM, 1992.

\bibitem[Hol95]{HOLSBOOKWAVE}
M.~Holschneider.
\newblock {\em Wavelets: An Analysis Tool}.
\newblock Oxford University Press, 1995.

\bibitem[H{\"o}r82]{HORMANDERI}
L.~H{\"o}rmander.
\newblock {\em The Analysis of Partial Differential Operators I}.
\newblock Springer Verlag, 1982.

\bibitem[Mey90]{MEYERIAII}
Y.~Meyer.
\newblock {\em Ondelettes et Op\`erateurs I\&II}.
\newblock Hermann, 1990.

\bibitem[Tri84]{TRIEBELII}
H.~Triebel.
\newblock {\em Theory of Function Spaces II}.
\newblock Birkhauser Verlag, 1984.

\end{thebibliography}
\end{document}